\documentclass[12pt,oneside]{amsart} 
\usepackage{amsfonts,newlfont,latexsym}
\usepackage{euscript}
\usepackage{amssymb,amscd}

      \theoremstyle{plain}
      \newtheorem{theorem}{Theorem}[section]
      \newtheorem{lemma}[theorem]{Lemma}
      \newtheorem{corollary}[theorem]{Corollary}
      \newtheorem{proposition}[theorem]{Proposition}
      \newtheorem{remark}[theorem]{Remark}
      
      \newtheorem{definition}[theorem]{Definition}      

      \makeatletter
      \def\@setcopyright{}
      \def\serieslogo@{}
      \makeatother

\def \R{\mathbb R}
\def \Rk{\mathbb R^k}

\def \Tk{\mathbb T^k}
\def \Tn{\mathbb T^n}
\def \Z{\mathbb Z}
\def \Zk{\mathbb Z^k}

\def \Es{E^{s}}

\def \Ci{C^{\infty}}

\def \a{\alpha}
\def \ta{\tilde \alpha}
\def \b{\beta}

\def \n{\mathfrak{n}}
\def \E{\mathfrak{e}}
\def \e{\varepsilon}
\def \A{\cal A}
\def \B{\cal B}
\def \H{\cal H}

\def \M{\cal M}
\def \S{\tilde {\cal M}}
\def \N{\cal N}
\def \L{\cal L}

\def \oa{{\cal O ^\ast}}
\def \o{{\cal O}}

\def\dist{\text{dist}}
\def \dim{\text{dim}\,}
\def \ker{\text{ker\,}}

\def\proof{{\bf Proof. }}
\def\QED{$\hfill\hfill{\square}$}

\numberwithin{equation}{section}

\begin{document}

\author{Boris Kalinin and Victoria Sadovskaya$^\ast$}

\address{Department of Mathematics \& Statistics,
 University of South  Alabama,   Mobile, AL 36688, USA}
\email{kalinin@jaguar1.usouthal.edu, sadovska@jaguar1.usouthal.edu}

\title [On the classification of resonance-free Anosov $\Z^k$ actions]
{On the classification of resonance-free Anosov $\Z^k$ actions} 

\thanks{$^\ast$ Supported in part by NSF grant DMS-0401014}

\begin{abstract} 
We consider actions of $\Zk$, $k \ge 2$, by Anosov diffeomorphisms 
which are uniformly quasiconformal on each coarse Lyapunov distribution. 
These actions generalize Cartan actions for which 
coarse Lyapunov distributions are one-dimensional. 
We show that, under certain non-resonance assumptions on the Lyapunov 
exponents, a finite cover of such an action is smoothly conjugate to an action 
by toral automorphisms.

\end{abstract}

\maketitle 


\section{Introduction}

In this paper we consider the problem of classification of higher-rank abelian
Anosov actions. These are actions of $\Zk$ and $\Rk$, $k\ge 2$, by commuting 
diffeomorphisms of a compact manifold with at least one Anosov element.
It is conjectured that all  irreducible actions of this type are {\em smoothly} 
conjugate to algebraic models \cite{KaSp}. 
We will concentrate on $\Zk$ actions.  In this case, the algebraic models are 
actions by commuting  automorphisms of tori and infranilmanifolds.
\vskip.1cm

We note the difference with the case of  Anosov actions of $\Z$, i.e. Anosov 
diffeomorphisms. A well-known conjecture states that they are 
{\em topologically} conjugate to automorphisms of tori and infranilmanifolds.
The conjugacy is typically only H\"older continuous, and thus a smooth
classification of  Anosov diffeomorphisms is impossible. Even the topological
classification was obtained only for diffeomorphisms of infranilmanifolds and for 
diffeomorphisms  of codimension one \cite{Fr, M, N}.
\vskip.1cm

The study of  higher-rank abelian actions originated in the
work on Zimmer's conjecture that the standard action of $SL(n,\Z)$ on $\Tn$, 
$n \ge 3$, is locally rigid. This means that  any $C^1$ small perturbation is 
smoothly conjugate to the original action. In the proof of the conjecture, 
the smoothness of the conjugacy was established using the action of a 
diagonalizable subgroup isomorphic to $\Z ^{n-1}$ \cite{H1,KL2}. 
This motivated Katok and Spatzier to study local rigidity for abelian actions.
In \cite{KS97} they established local rigidity for a broad class of {\em algebraic}
Anosov actions of $\Zk$ and $\Rk$, $k\ge 2$. The natural next step is to obtain
a smooth conjugacy to an algebraic model for actions more general than $C^1$ 
small perturbations of these models. This is referred to as
global rigidity or classification.

In  \cite{KL2} Katok and Lewis established the first global rigidity result for certain
actions of $\Z ^{n-1}$ on $\Tn$. The actions they considered have $n$ transversal 
one-dimensional foliations which are stable foliations for some Anosov elements. 
Such actions are called maximal Cartan. Recently, this result was extended by  
Rodriguez Hertz \cite{RH} to a broader class of maximal actions on tori assuming 
the existence of one Anosov element. We note that the case of a torus 
(or infranilmanifold) is special. In this case, an Anosov element is known to be 
topologically conjugate to an automorphism, which gives topological conjugacy of 
the whole action to an algebraic model. Thus the problem  reduces to showing the 
smoothness of this conjugacy.

For actions on an arbitrary manifold, one does not immediately have a conjugacy
to an algebraic model, and even an algebraic model itself, which makes this case 
considerably more difficult. The algebraic structure is constructed using dynamical
objects such as coarse Lyapunov distributions. These are the finest non-trivial 
intersections of stable distributions of  various Anosov elements of the action. 
They naturally extend the notion of stable distributions for Anosov maps. 

A coarse Lyapunov distribution can be also be described in terms of Lyapunov 
exponents (see Section \ref{Lyapunov} for details). For a $\Z$ action generated 
by an Anosov diffeomorphism Lyapunov exponents given by the Multiplicative 
Ergodic Theorem can be viewed as linear functionals on $\Z$. In general, these
functionals depend on a given invariant measure, and the corresponding Lyapunov 
distributions are only measurable. Note, however, that all functionals on $\Z$ are proportional. Moreover, for {\em any} invariant measure the sum of all Lyapunov distributions for Lyapunov exponents proportional to a given negative exponent 
with positive constants of proportionality is precisely the stable distribution of the 
Anosov diffeomorphism. The Lyapunov exponents
for a higher rank abelian action are linear functionals on $\Rk$ or $\Zk$. A coarse 
Lyapunov distribution can be defined as a direct sum of Lyapunov distributions for all
Lyapunov exponents positively proportional to a given one.  To show that these 
distributions are continuous and independent of an invariant measure one needs 
to assume the existence of sufficiently many Anosov elements (see Proposition 
\ref{properties}). This assumption is present in all known classification results
when the manifold is not assumed to be a torus.  

Another property important in the study of higher rank abelian actions is TNS.
An action is called totally nonsymplectic (TNS) if there are no negatively
proportional Lyapunov functionals. This is equivalent to the fact that any pair of
coarse Lyapunov distributions is contracted by some element of the action. 
This assumption was introduced in \cite{KNT,NT}
for the study of cocycle rigidity and was also used in measure rigidity results
\cite{KS3,KaSp0}.  All discrete actions in global rigidity results \cite{KS2,KaSp,KL2} 
are TNS. In this paper the TNS assumption is used primarily to ensure that certain
elements act transitively. In particular, it yields sufficient irreducibility of the action.  
We note that without some irreducibility assumption there is no hope for any type
of rigidity, since one can take products of Anosov diffeomorphisms. 

The simplest higher rank actions are the ones with one-dimensional coarse 
Lyapunov distributions.
Such actions are called Cartan. They generalize the maximal Cartan actions in \cite{KL2}.
The first classification results on an arbitrary manifold were obtained in \cite{KaSp} 
for certain Cartan actions of $\Rk$ and $\Zk$ with $k \ge 3$. 
Our main result, Theorem \ref{main}, allows higher-dimentional coarse Lyapunov 
distributions under the condition of uniform quasiconformality. This means that
for any element $b$ of the action, all vectors in a given coarse Lyapunov distribution 
at a given point are expanded or contracted by the iterates of $b$ at essentially the 
same rate (see Section \ref{uqc}). This condition is 
trivially satisfied for Cartan actions. Compared to the result for $\Zk$ 
actions in \cite{KaSp}, we also reduce the rank requirement to the optimal $k \ge 2$ 
and weaken the assumptions on Lyapunov exponents. 

\begin{theorem} \label{main}
Let $\a$ be an  action of $\,\Zk$, $k \ge 2$, by  $\Ci$ diffeomorphisms of a compact 
connected smooth manifold $\M$. Suppose  that 
\begin{enumerate}
\item all non-trival elements of $\a$ are Anosov and at least one is transitive;

\item $\a$ is uniformly quasiconformal on each coarse Lyapunov distribution;

\item $\a$ is totally nonsymplectic (TNS);

\item for any Lyapunov functionals $\chi_i$, $\chi_j$, and $\chi_l$,
         the functional $\,(\chi_i - \chi_j)$ is not proportional to $\chi_l$.
 \end{enumerate}
Then a finite cover of $\a$ is $\Ci$ conjugate to 
a $\Zk$ action by affine automorphisms of a torus.
\end{theorem}

We note that under conditions (1), (2), and (3) of the theorem, the Lyapunov 
exponents do not depend on the measure, as follows from Theorem \ref{metric} below.
Thus it suffices to verify condition (4) for just one invariant measure.
This condition is equivalent to the following: 
\vskip.2cm

\hskip.6cm (4$'$) {\it For any Lyapunov functionals $\chi_i$, $\chi_j$, and $\chi_l$,}
$$\{ a\in\Rk :   \;  \chi_i (a)=\chi_j(a)\} \,\ne\, \ker \chi_l.$$
Clearly, this is  weaker that the ``non-degeneracy" assumption in  \cite{KaSp} that 
$\ker  \chi_i \cap \ker \chi_j$ is not contained in $\ker \chi_l$. We note that $\Zk$ 
actions by automorphisms of nonabelian infranilmanifolds always have resonances 
of the type $\chi _i = \chi_j + \chi_l$. Thus assumption (4), which is slightly stronger 
than the absence of such resonances, serves to exclude such actions.
\vskip.1cm

In \cite{KS2} we obtained rigidity results for similar actions under the assumption 
that the coarse Lyapunov foliations are pairwise jointly integrable. This  assumption
can be easily verified for actions on tori. Our approach to the proof of Theorem \ref{main} 
is different from the one in \cite{KS2}. First we construct a Riemannian metric with 
respect to which the contraction/expansion is given by the Lyapunov exponents. 
The following result extends Theorem 1.2 in \cite{KaSp} to the case
of higher dimensional coarse Lyapunov distributions. 

\begin{theorem} \label{metric}
Let $\a$ be a  smooth action of $\,\Zk$  on a compact connected smooth manifold $\M$ 
satisfying conditions (1), (2), and (3) of Theorem \ref{main}. Then the Lyapunov
functionals are the same for all invariant measures. Moreover, there exists 
a H\"older continuous Riemannian metric  on  $\M$ such that 
for any  $a \in \Zk$ and any Lyapunov functional $\chi$
\begin{equation}
\| Da(v)\|=e^{\chi (a)} \| v \|  \label{metric eq1}
\end{equation}
for any vector $v$ in the corresponding Lyapunov distribution. 
\end{theorem}

After the construction of the metric, our approach diverges from the one in 
\cite{KaSp}, which required rank $k \ge 3$.  We prove a special application of the 
$C^r$ Section Theorem, Proposition \ref{C^r application}, which is of
independent interest. We use this result to establish smoothness 
of coarse Lyapunov distributions. Then we prove smoothness of the metric
constructed in Theorem \ref{metric} and obtain a smooth  conjugacy to 
an algebraic model.
\vskip.2cm 

In Section \ref{preliminaries} we introduce the main notions and summarize 
the results that are used throughout this paper.
In Section \ref{proof metric} we prove Theorem \ref{metric}, and in Sections
\ref{C^r}, \ref{smooth metric}, and \ref{smooth conjugacy} we complete the proof 
of our main result, Theorem \ref{main}.


\section{Preliminaries}   \label{preliminaries}

Throughout the paper, the smoothness of diffeomorphisms, actions, 
and manifolds is assumed to be $\Ci$, even though all definitions and 
some of the results can be formulated  in lower regularity.

\subsection{Anosov actions of $\,\Zk$ and $\Rk$}

\label{Anosov actions} $\;$
\vskip.1cm 

Let $a$ be a diffeomorphism of a compact manifold $\M$. 
We recall that  $a$ is  {\em Anosov} if there exist a continuous $a$-invariant
decomposition  of the tangent bundle $T\M=E^s_a \oplus E^u_a$ and constants 
$K>0$, $\lambda>0$ such that for all $n\in \mathbb N$
 \begin{equation} \label{anosov}
   \begin{aligned}
  \| Da^n(v) \| \,\leq\, K e^{- \lambda n} \| v \|&
      \;\text{ for all }\,v \in E^s_a, \\
  \| Da^{-n}(v) \| \,\leq\, K e^{- \lambda n}\| v \|&
      \;\text{ for all }\,v \in E^u_a. 
  \end{aligned}
  \end{equation}
The distributions $E_a^s$ and $E_a^u$ are called the {\em stable} and {\em unstable} 
distributions of  $a$.
\vskip.2cm

Now we consider a $\Zk$ action $\alpha$ on a compact manifold $\M$
via diffeomorphisms.
The action is called {\em Anosov}$\,$ if there is an element which acts as an 
Anosov diffeomorphism. Throughout this paper, we will be using the same 
letter for an element of the acting group and for the corresponding diffeomorphism.

For a $\Zk$ action $\alpha$ there is an associated
$\Rk$ action $\ta$ on a manifold $\S$ given by the standard suspension 
construction \cite{KaK}. We will refer to $\ta$ as the {\em suspension} of $\a$. 
It generalizes the suspension flow of a diffeomorphism. Similarly,
the manifold $\S$ is a fibration over the ``time" torus $\Tk$ with the fiber $\M$.
If $\a$ is an Anosov $\Zk$ action then $\ta$ is an Anosov $\Rk$ action. An $\Rk$ 
action is called Anosov if some element $a \in \Rk$ is Anosov in the sense of the
following definition. 

\begin{definition} Let $\alpha$ be a smooth  action of $\,\Rk$
on a compact manifold $\M$. An element $a \in \Rk$ is called {\em Anosov} 
or {\em normally hyperbolic} for $\alpha$ if there exist positive constants $\lambda$, 
$K$ and a continuous $\alpha$-invariant splitting of the tangent bundle
$$
T\M = E^s _a \oplus E^u _a \oplus T\o  
$$
where $T\o$ is the tangent distribution of the $\Rk$-orbits, and \eqref{anosov} holds
for all $n\in \mathbb N$.
\end{definition}

Both in the discrete and the continuous case it is well-known that the distributions 
$E_a^s$ and $E_a^u$  are H\"{o}lder continuous and tangent 
to the stable and unstable foliations $W_a^s$ and $W_a^u$ respectively \cite{HPS}.
The leaves of these foliations are $C^\infty$ injectively immersed Euclidean 
spaces. Locally, the immersions vary continuously in $\Ci$ topology. Such 
foliations are said to have {\em uniformly $\Ci$ leaves}. In general, the 
distributions $E^s$ and $E^u$ are only H\"older continuous transversally 
to the corresponding foliations.

The set of Anosov elements ${\cal A}$ in $\Rk$ is always an open subset 
of $\Rk$ by the Structural Stability Theorem for normally hyperbolic maps 
by Hirsch, Pugh and Shub \cite{HPS}. An $\Rk$ action is called 
{\em totally Anosov} if the set of Anosov elements $\cal A$ is dense in $\Rk$. 
Algebraic Anosov actions are always totally Anosov, however, this is not
known for nonalgebraic actions.

\subsection{Lyapunov exponents and coarse Lyapunov distributions}
\label{Lyapunov} $\;$
\vskip.1cm 

 We will concentrate on the case of $\Rk$ actions, the case of $\Zk$ is similar.
 We refer to \cite{KaSp} and \cite{KS2} for more details. 

Let  $a$ be a diffeomorphism of a compact manifold $\M$ preserving an 
ergodic  probability measure $\mu$. By Oseledec Multiplicative Ergodic 
Theorem, there exist finitely many numbers $\chi _i $ and a measurable 
splitting of the tangent bundle $T\M = \bigoplus E_i$ on a set of full measure 
such that the forward and backward Lyapunov exponents of $v \in  E_i $ 
are $\chi _i $.  This splitting is called {\em Lyapunov decomposition}.

Let $\mu$ be an ergodic probability measure for an $\Rk$ action $\alpha$
on a compact manifold $\M$. By commutativity, the  Lyapunov decompositions 
for  individual elements of $\Rk$ can be refined to a joint invariant splitting for 
the action. The following proposition from \cite{KaSp} describes the 
Multiplicative Ergodic Theorem for this case. See  \cite{KS2} for the discrete
time version and \cite{KaK} for more details on the Multiplicative Ergodic 
Theorem and related notions for higher rank abelian actions. 

 \begin{proposition} Let $\a$ be a smooth action of  $\Rk$ and let $\mu$ be an 
 ergodic invariant  measure.
 There are finitely many linear functionals $\chi $ on $\Rk$, a set of
 full measure ${\cal P}$, and an $\alpha$-invariant  measurable splitting 
 of the tangent bundle  $T\M = \bigoplus E_{\chi}$ over ${\cal P}$  such that 
  for all $a \in \Rk$ and $ v \in E_{\chi}$, the Lyapunov exponent of $v$ is $\chi 
  (a)$, i.e.
 $$
   \lim _{n \rightarrow \stackrel{+}{_{-}} \infty } 
   \frac{1}{n} \log \| D a ^n  (v) \| = \chi (a),
 $$ 
  where $\| .. \|$ is a continuous norm on $T\M$. 
\end{proposition}
 
The splitting  $ \bigoplus E_{\chi} $ is called the {\em Lyapunov decomposition}, 
and the linear functionals $\chi$ are called the {\em Lyapunov exponents} or 
{\em Lyapunov functionals} of  $\alpha$. The hyperplanes $\,\ker \,\chi \subset \Rk$ 
are called the {\em Lyapunov hyperplanes} or {\em Weyl chamber walls}, and the 
connected components of 
$\,\Rk - \cup _{\chi} \ker \chi $ are called the {\em Weyl chambers} of $\alpha$. 
The elements in the union of the Lyapunov hyperplanes are called 
{\em singular}, and the elements in the union of the Weyl chambers 
are called {\em regular}.
We note that the corresponding notions for a $\Zk$ action and for its suspension are 
directly related. In particular, the nontrivial Lyapunov exponents are the same.
In addition, for the suspension there is one identically zero Lyapunov exponent 
corresponding to the orbit distribution. From now on, the term Lyapunov exponent 
will always refer to the nontrivial functionals.

Consider  a $\Zk$ action by {\em automorphisms} of a torus or an infranilmanifold.
In this case, the Lyapunov decomposition is determined by the eigenspaces of the 
automorphisms, and the Lyapunov exponents are the logarithms of the moduli 
of the eigenvalues. Hence they are independent of the invariant measure, and 
they give uniform estimates of expansion and contraction rates. Also, every Lyapunov
distribution is smooth and integrable. 

In the nonalgebraic case, the individual Lyapunov distributions are in general 
only measurable and depend on the given measure. This can be observed already 
for a single Anosov diffeomorphism. However, the full stable distribution for any 
measure always agrees with $E^s_a$. For higher rank actions, {\em coarse Lyapunov distributions} play a similar role.
For any Lyapunov functional $\chi$ the coarse Lyapunov distribution is the direct 
sum of all Lyapunov spaces with Lyapunov functionals positively proportional to $\chi$:
 $$
    E^{\chi} = \oplus E_{\chi '}, \quad \chi ' = c \, \chi \,\text{ with }\, c>0. 
 $$   

One can see that for an algebraic action such a distribution is a finest nontrivial 
intersection of the stable distributions of certain Anosov elements of the action.
For nonalgebraic actions, however, it is not a priori clear that the intersection 
of several stable distributions has constant dimension and that it is better than
measurable. It was shown in \cite{KaSp} that, in the presence of sufficiently many 
Anosov elements,  the coarse Lyapunov distributions are indeed well-behaved
dynamical objects. The next proposition summarizes important results  
 for the suspensions of $\Zk$ actions under consideration
(see Proposition 2.4, Lemma 2.5, Corollary 2.8 in \cite{KaSp}).

\begin{proposition} \label{properties}
Let $\alpha$ be a $\Ci$ action of $\,\Zk$  on $\M$ for which all non-trival elements are 
Anosov and at least one is transitive. Let $\ta$ be its suspension $\Rk$ action 
on the manifold $\S$. Then 

\begin{enumerate}
\item There is  $\ta$-invariant H\"{o}lder continuous coarse 
Lyapunov splitting
$$
T\S=  T\o  \oplus \bigoplus E^i \qquad (\ast)
$$
where $T\o$ is the tangent distribution of the $\Rk$-orbits and 
$E^i$ are the finest nontrivial intersections of the stable distributions 
of various  Anosov elements of the action. Every  distribution $E^i$ is 
tangent to the H\"older foliation $W^i$ with uniformly $\Ci$ leaves.

\item The Lyapunov hyperplanes and Weyl chambers are the same
for all $\ta$-invariant ergodic probability measures. The set $\A$ of 
Anosov elements for $\ta$ is the union of these Weyl chambers in $\,\Rk$.
In particular, $\ta$ is a totally Anosov $\,\Rk$ action.

\item For any $\ta$-invariant ergodic probability measure the coarse 
Lyapunov splitting coincides on the set of full measure with the splitting 
$(\ast)$. More precisely, for each Lyapunov exponent $\chi$ 
 $$
   E^{\chi}(p) =  \bigcap _{\{a \in \A\, | \; \chi (a) <0\}} \Es _a(p) \;= 
   \bigoplus_ {\{ \chi ' = c \, \chi\,|\;  c>0 \}}  E_{\chi '} (p) 
 $$ 
 
\item The action is called {\bf totally nonsymplectic}, or {\bf TNS}, if there are 
no negatively proportional Lyapunov functionals. Under this condition,
almost every element of every Lyapunov hyperplane is transitive on $\S$.

\end{enumerate}

\end{proposition}

Property (2) implies that any regular element $a\in \R^k$ is Anosov. 
It is clear that its (un)stable distribution, given by Definition 2.1,
is the sum of all coarse Lyapunov distributions with $\chi(a)<0\;$ ($\chi(a)>0$).
For a singular element $a$ of  $\ta$,  we can define its neutral, stable, and 
unstable distributions  as follows:
$$
E_a^0 =T\o \oplus \bigoplus_ {\chi (a)=0} E^{\chi}, \qquad
E^s _a = \bigoplus_ {\chi (a) <0} E^{\chi}, \qquad 
E^u _a = \bigoplus_ {\chi (a) >0} E^{\chi}.
$$
Note that, a priori, we do not have any uniform estimates on the possible expansion 
or contraction of $E_a^0$ by $a$, so we cannot say that $a$ is  a partially hyperbolic 
element in the usual sense. However, the  properties of the action imply the following.

\begin{lemma} (\cite{KaSp}, Lemma 2.6) \label{SingularElements}
For the action $\ta$ in Proposition \ref{properties} the distributions $E_a^0$, $E^s _a$, 
and $E^u _a$ are H\"older continuous. $E^s _a$ and $E^u _a$ integrate to H\"older 
continuous foliations $W^s_a$ and $W^u_a$ with uniformly $\Ci$ leaves. $E^s _a$ 
is uniformly contracted and $E^u _a$ is uniformly expanded by $a$.
\end{lemma}

\subsection{Uniform quasiconformality, conformality, and conformal structures} \label{uqc}
$\;$ \vskip.1cm 

Let $a$ be a diffeomorphism of a compact Riemannian manifold 
$\M$, and let $E$ be a continuous $a$-invariant distribution. 
The diffeomorphism $a$ is {\em uniformly quasiconformal} on $E$ 
if  there exists a constant $K$ such that for all $n\in\Z$ and $x\in \M$
\begin{equation}\label{QC}
  K^E(x,a^n)=\frac{\max\,\{\,\|\,Da^n(v)\,\|\, :\; v\in E(x), \;\|v\|=1\,\}}
            {\,\min\,\{\,\|\,Da^n(v)\,\|\, :\; v\in E(x), \;\|v\|=1\,\}} \,\le\, K.
\end{equation}            
If $K^E(x,a)=1$ for all $x$, the diffeomorphism is said to be 
{\em conformal} on $E$.

We note that the notion of uniform quasiconformality does not depend on 
the choice of a Riemannian metric on the manifold.
Clearly, an Anosov diffeomorphism can be uniformly quasiconformal 
on $E$ only if $E$ is contained in its stable or its unstable 
distribution.

\vskip.2cm

A {\em conformal structure} on $\R^n$, $n\geq 2$, is a class of proportional 
inner products. The space $C^n$ of conformal structures on $\R^n$
identifies with the space of real symmetric 
positive definite $n\times n$ matrices with determinant 1, which is
isomorphic to $SL(n,\R) /SO(n,\R)$. It is known that the space 
$C^n=SL(n,\R) /SO(n,\R)$ carries a $GL(n,\R)$-invariant metric for 
which $C^n$ is a Riemannian symmetric space of non-positive curvature.
(See \cite{T} for more details.)

Let $E$ be a distribution on a compact manifold $\M$. 
For each $x\in \M$, let $C^E(x)$ be the space of conformal 
structures on $E(x)$. This gives a bundle $C^E$ 
over $\M$ whose fiber over $x$ is $C^E(x)$. 
A continuous (smooth, measurable) section of $C^E$ is called a continuous 
(smooth, measurable) conformal structure on $E$.
A measurable conformal structure $\tau$ on $E$ is called  bounded
if the distance between $\tau(x)$ and $\tau_0(x)$ is uniformly 
bounded on $\M$ for some continuous conformal structure $\tau_0$ 
on $E$.   

Clearly, a diffeomorphism is conformal with respect to 
a Riemannian metric on $E$ if and only if it preserves the conformal
structure associated with this metric.
\vskip.2cm

Let $\Gamma$ be a group acting on $\M$ via diffeomorphisms,
and let $E$ be a continuous $\Gamma$-invariant distribution.
We say that the action is uniformly quasiconformal on $E$
if the quasiconformal distortion $K^E(x, \gamma)$ is uniformly 
bounded for all $x \in \M$ and $\gamma \in \Gamma$. The next proposition
easily follows from the definition.

\begin{proposition} (\cite{KS2} Proposition 2.8) Suppose that the $\Gamma$-action 
 is generated by 
finitely many commuting diffeomorphisms. If each generator is uniformly 
quasiconformal then the  $\Gamma$-action is uniformly quasiconformal.

\end{proposition}

\begin{remark} \label{not proportional} 
If a $\Zk$ or an $\Rk$ action  is uniformly quasiconformal on a coarse 
Lyapunov distribution, then there is only one Lyapunov exponent
corresponding to this distribution for any given ergodic measure. 
Therefore, for the actions in Theorem \ref{main} and \ref{metric}
there are no positively proportional Lyapunov exponents,
i.e. coarse Lyapunov distributions coincide with Lyapunov distributions. 
Together with the TNS assumption, this also means that there are no proportional 
Lyapunov exponents.

\end{remark}

We note that uniform quasiconformality on a given coarse Lyapunov distribution
is stronger than the existence of only one Lyapunov exponent. This can be
observed even for algebraic actions with nontrivial Jordan blocks.


\section{Proofs}   \label{proofs}

\subsection{Outline of the proof} \label{outline} $\;$
\vskip.1cm 

In  Section \ref{proof metric} we obtain a Riemannian metric with respect to 
which the contraction/expansion is given by the Lyapunov exponents. The metric 
is constructed separately on each coarse Lyapunov distribution $E$. First we 
obtain an invariant conformal structure on $E$. Then we choose a proper 
normalization using a Livsic type argument for a transitive element in the Lyapunov 
hyperplane corresponding to $E$. 
This argument is specific to TNS case and is simpler then the one in \cite{KaSp}.

In Section \ref{C^r} we prove  a special  application of the $C^r$ Section
 Theorem of Hirsch, Pugh, and Shub. We apply
this result with one Lyapunov foliation as a base and a sum of several 
Lyapunov distributions as a fiber. This allows us to prove that any Lyapunov
distribution is smooth along any Lyapunov foliation and thus obtain the
smoothness of the Lyapunov splitting. 

In Section \ref{smooth metric} we establish smoothness of the metric
constructed in Theorem \ref{metric}. For each coarse Lyapunov distribution,
we use holonomy maps to show  smoothness of the conformal structure, 
and Livsic Theorem to  obtain  smoothness of the normalization.

 In Section \ref{smooth conjugacy} we show
that the metric gives rise to a smooth affine connection invariant under the action. 
This yields smooth conjugacy to an action by automorphisms of an infranilmanifold.
 We complete the proof by showing that this infranilmanifold is finitely covered by 
 a torus.


\subsection{Proof of Theorem \ref{metric}} \label{proof metric}$\;$
\vskip.1cm 
In this section we prove the following theorem which gives the corresponding 
Riemannian metric on the suspension manifold $\S$. Clearly, such a metric  
induces the desired metric on the original manifold $\M$ and implies Theorem \ref{metric}.

\begin{theorem} \label{metricS}
Let $\a$ be a  smooth action  of $\,\Zk$ on a compact connected smooth manifold $\M$ 
satisfying conditions (1), (2), and (3) of Theorem \ref{main}.
Then the Lyapunov exponents are the same for all invariant measures. 
Moreover, there exists a H\"older continuous Riemannian metric  on the 
suspension manifold $\S$ such that for any $a \in \Rk$ 
 and any Lyapunov exponent $\chi$
\begin{equation}
\| Da(v)\|=e^{\chi (a)} \| v\|  \label{metric eq}
\end{equation}
for any vector $v$ in the corresponding Lyapunov distribution. 
\end{theorem}

First we observe that the orbit distribution carries a natural metric invariant under 
the action. It is induced by the local diffeomorphism between an orbit and $\Rk$.
We can also declare that the coarse
Lyapunov distributions and the orbit distribution are pairwise orthogonal. Thus the 
theorem is reduced to constructing the desired metric separately on each coarse 
Lyapunov distribution. 

For the rest of the proof of the theorem, we fix a coarse Lyapunov distribution
$E$ and the  corresponding coarse Lyapunov foliation $W$.
We will construct the  Riemannian metric on $E$ satisfying the 
condition \eqref{metric eq}.

The first step is the construction of a conformal structure $\tau$ on $E$ invariant 
under $\ta$. The following theorem gives the corresponding
structure for the original action $\a$. We note that if $E$ is one-dimensional this
step becomes trivial.

\begin{theorem} \label{conformal} 
Let $\a$ be $\Ci$  action of $\Zk$ on a compact connected smooth manifold $\M$  
with a transitive Anosov element. Suppose that the action is uniformly 
quasiconformal on an invariant H\"older continuous distribution $E$.

Then  the action preserves a conformal structure $\tau$ on $E$ which is H\"older 
continuous on $\,\M$. If, in addition, the distribution $E$ is tangential to a 
foliation $W$ with uniformly $\Ci$ leaves, then $\tau$ is uniformly 
$\Ci$ along the leaves of $W$. 
\end{theorem}

We say that a function is {\em uniformly $C^r$} along the leaves of $W$
if its derivatives up to order $r$ in the directions of $W$ exist and
are continuous on the manifold.
\vskip.2cm

\proof
The proof of this theorem is similar  to the proof of Theorem 1.3 in \cite{S}.
There the corresponding result is obtained for the case of a single Anosov 
diffeomorphism (or flow) with $E$ being the (strong) stable distribution.
We outline the main steps of the proof for $\Zk$ actions.

First we use an argument by Tukia \cite{T} to obtain an invariant bounded 
measurable conformal structure on $E$. We outline the proof for our case.

\begin{lemma}  Let $\Gamma$ be a group acting on 
a Riemanian manifold $\M$ via diffeomorphisms. 
Suppose that  the action is uniformly quasiconformal 
on a continuous invariant distribution $E$. 
Then there exists a $\Gamma$-invariant bounded measurable conformal 
structure $\tilde \tau$ on $E$.
\end{lemma}

\proof 
Let $\tau_0$ be a continuous conformal structure on $E$. 
For each point $x$ in $\M$, we consider the set 
$S(x)=\{\,\gamma ^{-1}_*\, \tau_0(\gamma x), \;\gamma \in \Gamma\,\}$ 
of pull-backs of conformal structures on the orbit.
Since the action is uniformly quasiconformal,  $S(x)$ is a bounded 
subset of $C^E(x)$, the space of conformal structures on $E(x)$.
Since $C^E(x)$ has non-positive curvature, there exists a unique ball of 
the smallest radius containing $S(x)$.
We denote its center by $\tilde \tau (x)$. One can show that the conformal 
structure $\tilde \tau$ is invariant, bounded, and measurable. 
\QED
\vskip.3cm

The next proposition shows that this invariant conformal structure
can be made H\"older continuous by changing it on a set of measure zero.
Its proof is virtually identical to the proof of 
Proposition 3.2 in \cite{S}. 

\begin{proposition} \label{Holder structure}
Let $\tilde \tau$ be a  bounded measurable conformal structure on 
a H\"older continuous distribution $E$.  Suppose that $\tilde \tau$ is preserved 
by a transitive Anosov diffeomorphism $a$.
Then there exists a H\"older continuous $a$-invariant conformal 
 structure $\tau$ on $E$, which 
 coincides with $\tilde \tau$ on a set of full Bowen-Margulis measure.
\end{proposition}

We recall that the Bowen-Margulis measure $\mu$ for $a$ is the unique measure 
of maximal entropy. Since any other element  $b\in \Zk$ commutes with $a$,
the measure $b_* \mu$ is also $a$-invariant and has the same entropy as $\mu$.
Hence $b_* \mu=\mu$, and $\mu$ is invariant under the action.
This implies that we can choose a set of full  measure invariant under the whole 
action where $\tau$ coincides with $\tilde \tau$.
Since any such set is dense, we conclude that $\tau$ is $\a$-invariant.
\vskip.2cm

Now suppose that $E$ is tangential to a foliation $W$ with uniformly $\Ci$ leaves.
As in the proof of Theorem 1.3 in \cite{S},  one can show the smoothness of 
the conformal structure $\tau$ along the leaves of  $W$ as follows. 
Let $b$ be an element of the action which contracts $W$.
We recall that there exists a family of diffeomorphisms $h_x:W(x) \to E(x)$ 
which give a non-stationary linearization of $b$ along $W$ 
(see Proposition 4.1 in \cite{S}).
One can show that $h_x$ maps the restriction of $\tau$ to a leaf $W(x)$  
to the conformal structure $\tau(x)$ on $E(x)$ (see Lemma 3.1 in \cite{S}),
and hence $\tau$  is $\Ci$ along the leaves of  $W$.
This completes the proof of Theorem \ref{conformal}.

\QED

Theorem \ref{conformal} gives an invariant conformal structure for the 
original $\Zk$ action $\a$. This conformal structure gives rise to the conformal 
structure $\tau$ on the coarse Lyapunov distribution $E$ for the suspension 
$\tilde \alpha$. We fix a normalization of $\tau$ by a smooth positive function 
to obtain a background metric $g_0$ on $E$ with respect to which  $\tilde \alpha$ 
is conformal. The rest of the proof of Theorem \ref{metricS} consists of constructing 
a new metric $g$ in the conformal class of $g_0$ for which condition \eqref{metric eq} 
is satisfied. 
\vskip.1cm

For any element $b \in \Rk$ we denote by $D_x^E b$ the restriction of its derivative
at $x\in \S$ to $E(x)$. Let  $q(x,b)$ be the norm of  $D^E_x b$ induced by the metric 
$g_0$ on $E$. Since $E$ and $g_0$ are H\"older continuous, so is the function
$q(x,b)$  for any $b$. As the action is conformal
with respect to $g_0$, for any $b \in \Rk$ we have 
\begin{equation}\label{Holder0}
q(x,b) =  \| D^E_x b  \,\|_{g_0} =  \| (D_x b) (v) \|_{g_0}  \cdot \|v\|^{-1}_{g_0}
\;\text{ for any nonzero }v \in E(x).
\end{equation}
For a positive continuous function $\phi$ let $g$ be the new metric $E$ such that 
$\|v\|_g=\phi(x) \|v\|_{g_0}$ for any $v\in E(x)$. We will write $g=\phi g_0$. Then
$$
q(x,b) = \phi(x)\, \phi(bx)^{-1} \| (D_x b) (v) \|_g \cdot \|v\|^{-1}_g,
$$
and the condition \eqref{metric eq} is satisfied 
 for an element $b\in \Rk$  with respect to $g$ if and only if
\begin{equation}\label{Holder0.5}
   \phi (x) \cdot \phi (bx)^{-1} = e^{-\chi(b)} q(x,b) \;\text{ for any }  x\in \tilde \M.
\end{equation}

\begin{definition} An element $a \in \Rk$ is called {\em generic singular} if it is
contained in exactly one Lyapunov hyperplane.
\end{definition}

We take a generic singular element $a_0$ in the Lyapunov hyperplane $\L$ 
corresponding to $E$. 
We will construct a metric $g$ which satisfies condition \eqref{metric eq}, or 
equivalently \eqref{Holder0.5}, for the one-parameter subgroup $\{ta_0 \}$.
Any element $a$ of this subgroup  is neutral on $E$, i.e. $\chi (a)=0$.
Therefore, condition \eqref{metric eq} simply means that $g$ is preserved by this 
subgroup. 
By Proposition \ref{properties} (4), we can choose $a_0$ to be transitive on $\tilde \M$, 
so there exists a point $x^\ast$ with dense orbit $\oa =\{(ta_0)x^\ast \}$. 
We define a new metric $g^\ast= \phi g_0$ on $E$ 
over $\oa$ as the propagation of $g_0$ from $E(x^\ast)$ along this orbit by the 
derivative of  element $t a_0$, i.e. we choose $\phi ((ta_0)x^\ast) = q(x^\ast, ta_0)^{-1}$.
By the construction, the metric $g^\ast$ is preserved by the one-parameter subgroup 
$\{ta_0 \}$. 

The main part of the proof is to show that $g^\ast$ is H\"older continuous on $\oa$. 
Then it extends to a H\"older continuous Riemannian metric $g=\phi g_0$ on $E$,
which  is also preserved by the subgroup $\{ta_0 \}$. 
Now we consider an arbitrary element $b \in  \Rk$. By commutativity, the push forward 
$b_\ast g$ of the metric $g$ is again preserved by $\{ta_0 \}$. Since $b$ is conformal
with respect to $g=\phi g_0$, we observe that $b_\ast g = \psi g$ for some positive
function $\psi$. Since $\oa$ is dense, it is easy to see that this function must be constant: 
$\psi(x) = \psi (x^\ast) =c$. We conclude that $b_\ast g = c \cdot g$ on $\S$, which implies 
that for any $x \in \S$ and any nonzero vector in $E(x)$ the Lyapunov exponent exists and 
equals $\log c$. This shows that the condition \eqref{metric eq} is satisfied for any $b \in \Rk$.
\medskip

To prove that the metric $g^\ast$ is H\"older continuous on $\oa$ it suffices to show 
that for any point $x \in \oa$ which returns close to itself under an element $a=ta_0$ 
the norm $q(x,a)$ is H\"older close to $1$. Let $\b >0$ be such that all coarse 
Lyapunov distributions are H\"older continuous  with exponent $\b$. We will show
that there exist positive constants $\e_0$ and $K$ such that if $\dist(x, ax) <\e_0$ 
for some $x\in \S$ and $a=ta_0$ with $t>1$ then 
\begin{equation}\label{Holder1}
|\log q(x, a)| < K \cdot \dist (x, ax) \, ^\b
\end{equation}

We will use the following lemma which can be viewed as a generalization of 
Anosov closing lemma to the case of a partially hyperbolic diffeomorphism 
with integrable neutral distribution. For a point with a close return it gives a nearby 
point which returns to its leaf of the neutral foliation. Note that the usual 
Anosov closing lemma is for a fully hyperbolic system so that the neutral 
foliation is trivial in the discrete time case (Anosov diffeomorphism) and is
the orbit foliation in the continuous time case (Anosov flow). The main difference
of the partially hyperbolic case is that, unlike the orbit foliation of a hyperbolic flow, the neutral foliation is not smooth in general. Therefore, we can not use the contracting 
mapping principle approach which relies on differentiability of the holonomy
of the neutral foliation. Instead, we use a topological fixed point theorem and 
deal with stable and unstable directions separately.

\begin{lemma} \label{closing}
Let $W$ be a coarse Lyapunov foliation, let $E=TW$, and let $a_0$ be a generic
singular element in the corresponding Lyapunov hyperplane. There exist positive 
constants $\e_0$, $C$, and $\lambda$
such that for any $x\in \S$ and $a=ta_0$ with  $\dist(x, ax) =\e <\e_0$, 
there exists $y\in \tilde M$ and $\delta \in \Rk$ such that 
\begin{enumerate}

\item $\dist \,(x,y) < C\e$,
\item $ (a+ \delta) y \in W(y) $,
\item $ \| \delta \| <C\e$,
\item $\dist\, ((sa_0)x, (sa_0)y) <C\e e^{- \lambda \min \{s,t-s\} }\, $ for any $\,0\le s \le t$,
\item $|\log q(x,a)- \log q(y,a)|< C\e^\b$, where $q(x,a)$ is given by \eqref{Holder0}.

\end{enumerate}
\end{lemma}

\proof 
Note that any element $a=ta_0$, $t\ne 0$, is also generic singular, and the 
neutral, stable, and unstable distributions $E^0_a$, $E^s_a$, and $E^u_a$ 
are the same for all $t > 0$. By Lemma \ref{SingularElements}, $E^u_a$
($E^s_a$) is uniformly expanded (contracted) by $a$ and integrates to the
foliation $W_a^u$ ($W_a^s$). The distribution $E^0_a$ is the direct sum
$T\o \oplus E$, and hence it integrates to foliation $\o W$. We note that 
$E^u_a \oplus E$ is also integrable as the unstable distribution of an Anosov
elements close to $a$. Therefore, $E^u_a \oplus E \oplus \o$ also integrates
to the foliation which we denote by $F^u$. 

First we find a point $z$ close to $x$ such that $az \in F^u(z)$. We consider 
the holonomy map of foliation $F^u$ between the stable leaves $W^s_a(x)$ 
and $W^s_a(ax)$. Let $B^s_r (x)$ denote the ball of radius $r$ in $W^s_a(x)$.
If $\e _0$  chosen sufficiently small then the balls $B^s_1(x)$ and $B^s_1(ax)$
are close and the holonomy map $H: B^s_1(ax) \to  W^s_a(x)$ along the 
leaves of $F^u$ is well-defined. Moreover, since the angles between the tangent 
spaces to $W^s$ and $F^u$ are uniformly bounded away from $0$, there exists
a constant $C_1$ such that $\dist (x, H(ax)) \le C_1 \e$  and $H$ is  close to an 
isometry: 
$$
| \dist (H(z'), H(z'')) -\dist (z',z'')| \le C_1 \e \quad \text{ for all } z',z'' \in B^s_1(ax).
$$
We can also assume that the time $t$ of the $\e$-return has to be large enough
so that $a$ contracts $W^s_a$ by at least a factor of 2. Now it is easy to see that 
for $B=B^s_{4C_1}$ we have $H (a (B)) \subset B$. Indeed, for any $z \in B$
$$
\begin{aligned}
&\dist (x,H(az)) \le \dist(x, H(ax)) + \dist(H(ax), H(az)) \le \\
&C_1 \e + \dist(ax,az) + C_1 \e
\le C_1 \e + \frac12 \dist(x,z) + C_1 \e \le 4C_1.
\end{aligned}
$$
Therefore, by Brouwer fixed point theorem, there exists a point $z \in B$ such that 
$H (a(z))=z$, i.e. $az \in F^u (z)$. 

Now we can apply a similar argument to find the desired point $y$. Since the leaves of
$\o W$ foliate the leaves of $F^u (z)$, we can consider the holonomy map $H$ 
of foliation $\o W$ from $W^u_a(z)$ to $W^u_a(az)$ inside the leaf $F^u (z)$. 
As above, we obtain a fixed point  for the map 
$H_1 \circ a^{-1}: B \subset W^u_a(az) \to W^u_a(z)$. This give the existence of
the point $y \in W^u_a(z)$ for which $ay \in \o W(y)$ with $\dist (z,y) \le C_2 \e$. Hence there
exists $\delta \in \Rk$ with $\| \delta \| <C_3 \e$ such that condition (2) is satisfied.

Since by the construction $z \in W_a^s(x)$ and $y\in W_a^u(z)$, one can easily
see that (4) is satisfied. Then (5) follows from the standard estimate for a H\"older
expansion/contraction coefficient along exponentially close trajectories. 
\QED

\vskip0.3cm

We will now complete the proof of the theorem by establishing \eqref{Holder1}.
We use the notations of the previous lemma and let $b=a+\delta$. Since $\delta$ 
is small, we have $|\log q(y,a)-\log q(y,b)| < C_1 \e$.
Together with Lemma \ref{closing} (5) this implies 
\begin{equation}\label{Holder2}
|\log q(x,a)-\log q(y,b)| < C_2 \e ^\b,
\end{equation}
which reduces the proof of \eqref{Holder1} to showing that
\begin{equation}\label{Holder3}
|\log q(y, b)| < C_3 \e ^\b
\end{equation}
To show this we take an Anosov element $c$ which contracts $W$. 
Let $y_*= \lim (t_n c) y$ be an accumulation point of the $c$-orbit of $y$.
We observe that $y_*$ is a fixed point for $b$. 
Indeed, using commutativity and the fact that $by \in W(y) \subset W^s_c (y)$
we obtain $by_*=\lim b((t_nc) y)=\lim (t_nc )(by)=y_*$.
Our goal is to show that $q(y_*, b)$ is close to $1$ and then to show that $q(y_*, b)$ 
is close to $q(y, b)$ to obtain \eqref{Holder3}. 

Since $\delta$ is small, we may assume that $b$ is not contained in any Lyapunov 
hyperplane, except possibly for the Lyapunov hyperplane $\L$ corresponding to $W$.
If $b \in \L$ then for the fixed point $y_*$ we must have $q(y_*,b)=1$. Indeed,
$q(y_*,b) \ne1$ is impossible since the Anosov elements close to $b$ on one side
of $\L$ contract $E$ while those on the other side of $\L$ expand $E$. If $b \notin \L$ 
then $b$ is regular and hence Anosov.  Then it is well known that the orbit $\Rk y_*$ of its 
fixed point $y_*$ must be compact \cite{Qian94}. 
The Lyapunov exponents and the Lyapunov splitting are defined everywhere on
this compact orbit. Let $\tilde \chi$ be the Lyapunov exponent on this orbit 
corresponding to $E$. By Proposition \ref{properties} we have 
$\ker \tilde \chi = \L$.  Thus $\tilde \chi (a)  =0$ and we obtain 
$$|\tilde \chi (b)| =|\tilde \chi (a) + \tilde \chi (\delta)|= |\tilde \chi (\delta)| < 
C_4 \|\delta\| < C_5  \e $$
We note that $\tilde \chi (b) = \log (q(y_*,b))$ since $y_*$ is a fixed point for $b$. 
Thus we conclude that in both cases, $b \in \L$ and $b \notin \L$, we have
\begin{equation}\label{Holder4}
|\log q(y_* , b)| < C_5 \e
\end{equation}
Now we will show that
 \begin{equation}\label{Holder5}
 |\log q(y,b)-\log q(y_* ,b)| \le C_6 \e^\b 
 \end{equation}
Using commutativity we can write $b=(-t_nc) \circ b \circ (t_nc)$ and obtain
$$q(y, b) = q(b(t_nc)y, -t_nc)   \cdot q ((t_nc)y, b)  \cdot q(y, t_nc)$$ 
Since $(t_nc)y \to y_*$ the middle term $q ((t_nc)y, b)$ tends to $q(y_*, b)$. 
On the other hand the product of the other two terms is close to 1:
$$
|\log \,[q(y, t_nc) \cdot q(b(t_nc)y, -t_nc) ]\, |=
| \log q(y, t_nc) - \log q(by, t_nc)|<C_6 \e^\b.
$$
For the last inequality we note that $c$-orbits of $y$ and $by$ are exponentially 
close and hence the contraction coefficients $q(by, t_nc)$ and $q(y, t_nc)$ are
H\"older close by the standard telescopic product argument. 
This completes the proof of \eqref{Holder5}. 

The inequalities \eqref{Holder5} and \eqref{Holder4} yeild \eqref{Holder3}, which together
with  \eqref{Holder2} gives the desired H\"older estimate \eqref{Holder1}. 
This completes the proof of the theorem.

\QED

\subsection{$C^r$ Section Theorem and smoothness of coarse Lyapunov splitting} 
\label {C^r} $\;$
\vskip.1cm 

In this section we prove that, under the assumptions of Theorem \ref{main},
all coarse Lyapunov distributions are $\Ci$. We will use the following a special 
version of  the $C^r$ Section Theorem of Hirsch, Pugh, and  Shub (see Theorems 
3.1, 3.2, and  3.5,  and Remarks 1 and 2 after Theorem 3.8 in  \cite{HPS}). 

\begin{theorem} \label{C^r section th}  (\cite{HPS}) Let $f$ be a $C^r$, $r\ge 1$, 
diffeomorphism of a compact smooth manifold $\M$. Let $W$ be an 
$f$-invariant topological foliation with uniformly $C^r$ leaves. 
Let $\B$ be a normed vector bundle over $\M$ and $F:\B \to \B$ be a linear 
extension of $f$ such that both $\B$ and $F$ are uniformly $C^r$ along 
the leaves of $W$.

Suppose that $F$ contracts fibers of $\B$, i.e. for any $x \in \M$ and any $u \in \B(x)$ 
\begin{equation} \label{k_x}
   \|F(u)\|_{fx} \le k_x \|u\|_x \;\text{ with }\, \sup_{x \in \M}  k_x <1.
\end{equation}
Then there  exists a unique continuous $F$-invariant  section of $\B$.
Moreover, if 
\begin{equation}
   \sup_{x \in \M}  k_x \a_x ^r <1 \;\text{ where }
    \a_x = \|(df|_{TW(x)})^{-1}\|, \label{C^r eq}
\end{equation}
then the unique invariant section is uniformly $C^r$ smooth along the leaves of $W$. 
\end{theorem}

\begin{remark} \label{C^r remark} Note that if $\a_x \le 1$ for all $x$, i.e. $f$ does 
not contract  $W$, then  \eqref{C^r eq} follows from \eqref{k_x}, and 
the invariant section is  uniformly $C^r$ along the leaves of $W$.
\end{remark}

In this theorem we consider the smoothness of the invariant section only
along the leaves of the invariant foliation $W$, so we need the 
smoothness of the bundle $\B$ and the extension $F$ only  along $W$. One can 
easily see from the proof in \cite{HPS} that the contraction in the base 
\eqref{C^r eq} needs to be estimated only along $W$. One  can also obtain 
this by formally applying the $C^r$ Section Theorem with the base manifold 
$\M$ considered as the disjoint union of the leaves of $W$. This has been observed 
in the study of partially hyperbolic systems (see the introduction, Theorem 3.1, 
and remarks after it in \cite{PSW}).

We use Theorem \ref{C^r section th} to obtain the following corollary on smoothness 
of invariant distributions. In contrast to the usual regularity results for Anosov and
partially hyperbolic systems, the distribution $E$ here is not $T\M$ and is not 
contained in $TW$.  This gives us greater flexibility that we will use to show
the smoothness of coarse Lyapunov distributions.

\begin{proposition}  \label{C^r application}  Let $f$ be a diffeomorphism 
of a compact smooth manifold $\M$. Let $W$ be an  $f$-invariant topological
foliation with uniformly $\Ci$ leaves   such that $\;\|(Df|_{TW(x)})^{-1}\| \le 1$ 
for all $x \in \M$. Let $E_1$ and $E_2$ be continuous 
$f$-invariant distributions on $\M$  such that the distribution $E=E_1 \oplus E_2$ 
is uniformly $\Ci$ along the leaves of $W$, and
for any $x$ in $\M$
  $$   
        \max  \,\{ \| Df(v)\| :\; v\in E_2(x), \;\|v\|=1 \} \, < \,
        \min   \,\{ \|Df(v)\| :\; v\in E_1(x), \;\|v\|=1 \}.
   $$
Then $E_1$ is uniformly $\Ci$  along the leaves of $W$.
\end{proposition}

\proof
We apply the $C^r$ Section Theorem \ref{C^r section th} and Remark
\ref{C^r remark}.
The distribution  $E$  is a continuous bundle
over $\M$ which is uniformly $\Ci$ along the leaves of $W$. There exist distributions
$\bar E_1$ and $\bar E_2$ in $E$ which are close to $E_1$ and $E_2$
respectively and  $\Ci$ along the leaves of $W$. Now we can consider a vector
bundle $\B$ whose fiber over $x$ is the set of linear operators from $\bar E_1(x)$ 
to  $\bar E_2(x)$ with the standard operator norm. The differential of $f$ induces a 
natural action $F$ on $\B$.   

 The assumptions on $f$ imply that $\sup_{x \in M}  k_x <1$, i.e. the induced map 
 $F$ contracts the fibers. 
Indeed, distributions $\bar E_1$ and $\bar E_2$ can be chosen sufficiently close to
$E_1$ and $E_2$, so that for all $x$ in $\M$
  $$   
        \max  \,\{ \| Df(v)\| :\; v\in \bar E_2(x), \;\|v\|=1 \} \, < \,
        \min   \,\{ \|Df(v)\| :\; v\in \bar E_1(x), \;\|v\|=1 \}.
   $$
Then for all $x$ in $\M$, we can take 
$$
      k_x= \frac {\max  \,\{ \| Df(v)\| :\; v\in \bar E_2(x), \;\|v\|=1 \}} 
       {\min   \,\{ \|Df(v)\| :\; v\in \bar E_1(x), \;\|v\|=1 \}} <1
 $$      
Since $k_x$ depends continuously on $x$, $\,\sup_{x \in M}  k_x <1$.
Thus there exists a continuous $F$-invariant section, and by uniqueness
the graphs of this section give the distribution $E_1$.
Moreover, since $\,\a_x =\|(df|_{TW(x)})^{-1}\| \le 1$, it follows that 
 the distribution $E^1$ is $C^\infty$ smooth along the leaves of $W$.
\QED
\vskip.3cm

Now we use Proposition \ref{C^r application} to show the smoothness of the 
coarse Lyapunov distributions. First we fix a coarse Lyapunov foliation $W$ 
and show that any coarse Lyapunov distribution  is smooth along its leaves.

\begin{proposition} \label{smooth distribution}  
Let $W$ be a coarse Lyapunov foliation. Under the assumptions of Theorem \ref{main},
any coarse Lyapunov distribution is uniformly $\Ci$ along the leaves of $W$. 
\end{proposition}

\proof
We consider a generic 2-dimensional subspace $\pi$ in $\Rk$ which intersects 
the Lyapunov hyperplanes along distinct lines.  For any Lyapunov exponent $\chi_j$
we denote by $\H_j$ the half-plane in $\pi$ on which  $\chi_j$ is negative.
We now order these half-planes and the corresponding Lyapunov
exponents counterclockwise so that $W$ is the Lyapunov foliation $W^1$ 
corresponding to $\chi_1$. 
Then there exists a unique $i>1$ such that $\H_1 \cap \H_j$ is non-empty
for $1 \le j \le i$, and $(-\H_1) \cap \H_j$ is non-empty for $j>i$. 
Note that by the TNS assumption (3),  $-\H_1\ne \H_j$ for any $j$.

We observe that $\,\bigoplus_{j=1}^i E^j=E^s_a\,$ for any 
Anosov element $a$ in $\H _1 \cap \H _i$. Therefore, this distribution is uniformly 
$\Ci$  along $W^s_a$, and in particular along $W^1$. We will apply 
Proposition \ref{C^r application} with $W=W^1$, $E=\bigoplus_{j=1}^i E^j$,
and various choices for $E_1$, $E_2$, and $f$. 

Recall that assumption (4) of the theorem is equivalent to (4$'$), and it
implies that the hyperplanes $\,\{ a\in\Rk :   \;  \chi_i (a)=\chi_j(a)\}\, $ and 
$\,\ker \chi_1\,$ do not coincide for distinct $ i$ and  $ j$. 
Thus, the plane $\pi$ can be chosen in such a way
that $\, \chi_i (b)\ne \chi_j(b)\,$ for any element $b\ne 0$  in $\,\ker \chi_1 \cap \pi\,$
and  any  $ i \ne  j$.
We take an element $b \in \ker \chi_1 \cap \pi$ such that $b\in \H_j$ for 
$ 1 < j \le i$.
Since the values $\chi _j (b)$, $1 \le j \le i$ are all different, we can reorder 
the indexes $1, ..., i$ so that 
$ \chi _{j_i} (b) < ... <\chi _{j_2} (b)< \chi _{j_1} (b)=0$, where $j_1=1$.
\vskip.1cm

We fix $m$, $1 \le m <i$, and apply Proposition \ref{C^r application}
with   $f=b$, $\,E_1=\,\bigoplus_{l=1}^m E^{j_l}$, and 
$E_2=\,\bigoplus_{l=m+1}^i E^{j_l}$. Note that  
    $$ \begin{aligned}
         \| Db(v)\| & \le e^{\chi_{m+1} (b)} \| v\| \;\text{ for any } v\in E_2, \text{ and} \\
         \| Db(v)\| & \ge e^{\chi_{m} (b)} \| v\| \;\text{ for any } v\in E_1
         \end{aligned}
     $$
with respect to the metric given by  Theorem \ref{metric}.
Since $\chi_{m+1} (b) < \chi_m (b)$,  Proposition \ref{C^r application} implies that 
$\,\bigoplus_{l=1}^m E^{j_l}$ is uniformly $\Ci$ along the leaves of $W^1$. 
\vskip.1cm

On the other hand, if we take $E_1=\bigoplus_{l=m+1}^i E^{j_l}$,
$\;E_2=\bigoplus_{l=1}^m E^{j_l}$, and  $f=b^{-1}$, we obtain that 
$\,\bigoplus_{l=m+1}^i E^{j_l}$ is uniformly $\Ci$ along  $W^1$. 
We conclude that for any $m$, $1 \le m <i$, both distributions 
$\,\bigoplus_{l=1}^m E^{j_l}$ and $\,\bigoplus_{l=m+1}^i E^{j_l}$
are uniformly $\Ci$ along  $W^1$. Since any coarse Lyapunov
distribution $E^j$, $1 \le j \le i$, can be obtained as an intersection 
of such distributions, it is also uniformly $\Ci$ along the leaves of $W^1$.

The smoothness of the distributions $E^j$, $j > i$, along $W^1$ can be obtained 
similarly.

\QED

\begin{corollary} \label{smooth splitting} Under the assumptions of 
Theorem \ref{main},  any coarse Lyapunov distribution is $\Ci$ on $\M$. 
In particular, the stable and unstable distributions of any element of the action
are $\Ci$.
\end{corollary}

\proof
Proposition \ref{smooth distribution}  implies that any coarse Lyapunov 
distribution $E$ is uniformly $\Ci$ along the leaves of any coarse Lyapunov 
foliation. To conclude that it is $\Ci$ on $\M$, we apply the following lemma
with $E$ considered as a function from $\M$ to the Grassmann bundle 
over $\M$ of subspaces of dimension $\dim E$.

\begin{lemma} \label{Ind Journe} 
Let $\phi$ be a map from $\M$ to a finite-dimensional $\Ci$ manifold.
If $\phi$ is uniformly $\Ci$ along the leaves of every coarse Lyapunov 
foliation then it is $\Ci$ on $\M$.
\end{lemma}

It is a well-known result  \cite{KL2, GS}
obtained by inductive application of Journ\'e Lemma \cite{J}.
\QED
 

\subsection{Smoothness of the metric constructed in Theorem \ref{metric}.} 
\label{smooth metric} $\;$
\vskip.1cm  
 
 \begin{proposition}
 Under the assumptions of Theorem \ref{main}, the metric $g$ 
 constructed in Theorem \ref{metric} is $\Ci$.
  \end{proposition}
 
 \begin{remark}
 As we already obtained the smoothness of the coarse Lyapunov splitting,
 we do not use assumption (4) in the proof of this proposition. 
 \end{remark}
 
 \proof
 We recall that we constructed the metric on each coarse Lyapunov 
 distribution separately, and then declared these distributions pairwise
 orthogonal. Since all coarse Lyapunov distributions are $\Ci$,
 it suffices to show that the metric on each coarse Lyapunov 
distribution is $\Ci$ on $\M$. First we will establish the smoothness of 
the conformal structure $\tau$ obtained in Theorem \ref{conformal}, and
then the smoothness of the normalization function $\phi$.

We fix a coarse Lyapunov distribution $E$ and the corresponding coarse 
Lyapunov foliation $W$. By Theorem \ref{conformal}, the conformal 
structure $\tau$ on $E$ is uniformly $\Ci$ along the leaves of $W$. We will 
now show that $\tau$ is $\Ci$ on $\M$. 

We consider Anosov elements $a,b \in \Zk$ such that $\chi (a) <0$, $\chi (b) >0$,
and they are not separated by any other Lyapunov hyperplane. We recall that under
the assumptions of the theorem no other
Lyapunov exponent has the same kernel (see Remark \ref{not proportional} ).
Thus for any other Lyapunov exponent $\chi '$, $\;\chi '(a)$ and $\chi' (b)$ have 
the same sign. It is easy to see that
$$
  E=E^s_a \cap E^u_b,  \quad\; E \oplus E^u_a = E^u_b,  \quad 
  \text{and} \quad E \oplus E^s_b = E^s_a.
$$

We will show  that $\tau$ is uniformly $\Ci$ along the leaves of $W^s_b$.
We note that $W^s_b$ is a subfoliation of $W^s_a$, and both foliations
 are $\Ci$ by Corollary \ref{smooth splitting}.
We consider two nearby points $x,y$ on the same leaf $W^s_b$ and the holonomy 
map $H_{x,y} : W(x) \to W(y)$ given by foliation $W^s_b$ within the leaf of $W^s_a$.
We will show that the map $H_{x,y}$ is conformal, i.e. it preserves the conformal 
structure $\tau$.   We give an argument similar to the one in the proof of Theorem 1.4 in \cite{S}. 
We can write
$$
H_{x,y} = a^{-n} \circ H_{a^nx,a^ny} \circ a^n,
$$
where $H_{a^nx,a^ny}$ is the corresponding holonomy from $W(a^nx)$ to $W(a^ny)$.
Since $\dist(a^nx,a^ny) \to 0$ as $n \to \infty$, the derivative $D_{a^nx} H_{a^nx,a^ny}$ 
of $H_{a^nx,a^ny}$ at $a^nx$ becomes close to an isometry. Since $a$ is conformal,
we conclude that the derivative of $a^{-n} \circ H_{a^nx,a^ny} \circ a^n$ at $x$ becomes 
close to conformal. Taking the limit, we obtain that the holonomy $H_{x,y}$ is conformal.
Since the foliation $W^s_b$ is smooth and $\tau$ is preserved by the corresponding
holonomy maps, we conclude that $\tau$ is uniformly $\Ci$ along  the leaves of $W^s_b$.

The smoothness along the leaves of $W^u_a$ is obtained similarly.
Since the conformal structure $\tau$ is  uniformly $\Ci$ along the leaves of $W^s_b$, 
$W^u_a$, and $W$, Lemma \ref{Ind Journe} yields the smoothness of $\tau$ on $\M$.
\vskip.1cm

To complete the proof of the proposition, it remains to show that 
the renormalization function $\phi$ from the proof of Theorem \ref{metric} is 
$\Ci$ on $\M$. We recall that for any $b \in \Zk$ 
and $x \in \M$ the function $\phi$ satisfies the cohomological equation
$$
\phi (x) \cdot \phi (bx) ^{-1} = e^{-\chi(b)}q(x,b)
$$ 
with $q(x,b)=\| Db|_{E(x)} \|_{g_0}$, where the metric $g_0$ is a
normalization of $\tau$. Since both $g_0$ and $E$ are $\Ci$,
the function $q(x,b)$ is $\Ci$ on $\M$. Now it follows from Livsic Theorem \cite{LMM} 
that the continuous solution $\phi$ is, in fact, $\Ci$.

\QED
 

 \subsection{Smooth conjugacy to a toral action} $\;$ 
 \label{smooth conjugacy}
\vskip.1cm

The smoothness of the coarse Lyapunov splitting and the existence of
smooth metric $g$ satisfying \eqref{metric eq1} give strong indication
of rigidity. If all coarse Lyapunov distributions are one-dimensional
the algebraic structure can be easily obtained since the vector fields of
unit length in coarse Lyapunov directions span a finite dimensional Lie
algebra \cite{KaSp}. In our case, it can be shown that the metric $g$ is 
flat. This, however, requires considerable effort. Instead, we use the metric 
$g$ to obtain $\a$-invariant  $\Ci$ affine connection on $\M$ and complete 
the proof  similarly to \cite{KS2}.

For each coarse Lyapunov foliation $W^i$ we denote by $\nabla^i$ 
the Levi-Civita connection of the induced Riemannian metric on the 
leaves of $W^i$. This connection is $\alpha$-invariant 
since any element of the action is a homotety on the leaves of $W^i$.
Since both the metric and $W^i$ are $\Ci$ on $\M$, so is the connection.
Now we define an affine connection $\nabla$ on $\M$ using a standard
construction. Let $X$ and $Y$ be two vector fields on $\M$. We use the 
coarse Lyapunov splitting $T\M=\oplus E^i$ to decompose 
$X=\sum X^i$ and $Y=\sum Y^i$, where $X^i, Y^i \in E^i$. Then
 $$
   \nabla_X Y = \sum_i  \nabla^i _{X^i} Y^i + \sum_{i \ne j} \Pi_j [X^i, Y^j],
 $$
where   $\Pi_j$ is the projection onto $E^j$, defines an 
affine connection. Since the distributions $E^i$ and the connections $\nabla^i$
are $\Ci$, $\nabla$ is also $\Ci$. Since $\nabla^i$ and $E^i$ are $\alpha$-invariant,
so is  $\nabla$. 
\vskip.1cm

Now we  consider a transitive Anosov element of  $\a$. 
This element preserves the  $\Ci$ affine connection $\nabla$ and has $C^\infty$   
stable and unstable distributions by Corollary \ref{smooth splitting}.
It is well known that such a diffeomorphism is conjugate to an Anosov 
automorphism of an infranilmanifold $\N$ by a $\Ci$ diffeomorphism $h$  \cite{BL}.  
Then $h$ conjugates the whole action $\a$ to an action $\bar \a$ by affine 
automorphisms of $\N$. This follows from the fact that any diffeomorphism 
commuting with an Anosov
automorphism of an  infranilmanifold is an affine automorphism itself (see \cite{H1}, 
proof of Proposition 2.18, and \cite{PY}, proof of Proposition 0). 

It remains to show that the infranilmanifold $\N$ is finitely covered by a torus. 
Recall that $\N$ is finitely covered by a nilmanifold $N/\Gamma$, where $N$ 
is a simply connected nilpotent Lie group, and $\Gamma$ is a cocompact lattice in $N$.
We need to show that $N$ is abelian. The Lie algebra $\n$ of $N$
splits into Lyapunov subspaces $\E_i$ with Lyapunov functionals $\chi_i$
of the action $\bar \a$. We note that since the conjugacy $h$
is smooth, the Lyapunov functionals of the actions $\a$ and $\bar \a$
coincide. In particular, $\bar \a$ has no proportional Lyapunov exponents
and no resonances of the type $\chi _l =\chi_i + \chi_j$ among the Lyapunov exponents.

If for $u \in \E_i$ and $v \in \E_j$ the bracket $[u,v] \ne 0$,
then $[u,v]$ belongs to a nontrivial Lyapunov subspace with Lyapunov functional
$\chi_i + \chi_j$. If $i=j$ this would give proportional Lyapunov functionals $\chi _i$
and $2\chi _i$, which is impossible. If $i \ne j$ this would give a resonance  
among the Lyapunov functionals $\chi _i$, $\chi _i$, and $\chi _l =\chi_i + \chi_j$, which 
is also impossible. 
Thus the Lie group $N$ is abelian, and the infranilmanifold $\N$ is finitely 
covered by a torus. This completes the proof of Theorem \ref{main}.


\end{document}